\documentclass{article}
\usepackage[utf8]{inputenc}
\usepackage{amsmath,amssymb,algorithm,algpseudocode,xcolor}
\usepackage{graphicx}
\usepackage{caption}
\usepackage{subcaption}
\usepackage[shortlabels]{enumitem}
\usepackage{hyperref}
\usepackage{cite}

\usepackage{tikz}
\usepackage{pgfplots}

\title{A new fast numerical method for the generalized Rosen-Zener model}

\author{C. Bonhomme\footnotemark[1], S. Pozza\footnotemark[2], N. Van Buggenhout \footnotemark[2]}
\date{}

\newtheorem{theorem}{Theorem}[section]
\newtheorem{lemma}[theorem]{Lemma}

\newlength\figureheight
\newlength\figurewidth

\usepackage{booktabs,subcaption,amsfonts,dcolumn}
\newcolumntype{d}[1]{D..{#1}}

\begin{document}

\maketitle

	\renewcommand{\thefootnote}{\fnsymbol{footnote}}
	\footnotetext[1]{Laboratoire Chimie de la Matière Condensée de Paris, LCMCP, Sorbonne Université, CNRS, 75005 Paris, France.}
\footnotetext[2]{Department of Numerical Mathematics, Charles University, Sokolovská 83, 186 75 Praha 8, Czech Republic.}
\footnotetext{This work was supported by Charles University Research programs UNCE/SCI/023 and PRIMUS/21/SCI/009 and by the Magica project ANR-20-CE29-0007 funded by the French National Research Agency.}

\begin{abstract}
   In quantum mechanics, the Rosen-Zener model represents a two-level quantum system. Its generalization to multiple degenerate sets of states leads to larger non-autonomous linear system of ordinary differential equations (ODEs).
   We propose a new method for computing the solution operator of this system of ODEs. 
   This new method is based on a recently introduced expression of the solution in terms of an infinite matrix equation, which can be efficiently approximated by combining truncation, fixed point iterations, and low-rank approximation.
   This expression is possible thanks to the so-called $\star$-product approach for linear ODEs. In the numerical experiments, the new method's computing time scales linearly with the model's size. We provide a first partial explanation of this linear behavior.
\end{abstract}

\section{Introduction}
Many crucial quantum mechanics problems require solving a time-dependent Schr\"odinger equation of the kind
\begin{equation}\label{eq:ODE}
    \frac{\partial}{\partial t} \psi(t) = -iH(t)\psi(t), \quad \psi(t_0) = \psi_0, \quad t \in I =[t_0,t_f],
\end{equation}
where the $N \times N$ matrix-valued function $H(t)$ -- the Hamiltonian -- is Hermitian and the solution $\psi(t)$ is the \emph{state vector} of the quantum system. When dealing with many-body problems, these systems are characterized by an exponential growth in size, i.e., the size of the system scales exponentially with the number of spins in, for instance, Nuclear Magnetic Resonance (NMR). As a consequence, computing the solution quickly becomes expensive in memory and computation cost – a problem known as the exponential-scaling wall \cite{Kup23}.
The problem becomes even more difficult when, instead of a single vector (state) solution, we want to solve the related equation for a matrix (operator) $U(t)$
\begin{equation}\label{eq:ODE:operator}
    \frac{\partial}{\partial t} U(t) = -iH(t)U(t), \quad U(t_0) = I_N, \quad t \in I =[t_0,t_f],
\end{equation}
where the unknown $U(t)$ is now an $N \times N$ matrix-valued function (the \emph{operator solution}) and $I_N$ is the identity matrix of size $N$. 
Naturally, the relation $\psi(t) = U(t) \psi_0$ connects the two systems in Equation \eqref{eq:ODE} and \eqref{eq:ODE:operator}.

Many algorithms have been proposed to solve such systems (we refer the reader to \cite{BlaCas17book,HaiHocIseLub06,IseMunNorZan00,Kop19}) and the research is very active in this field, pushing further the computational efficiency year by year \cite{BlaCasEsc22a,BlaCasChaEsc21,BerBlaCasEsc23,BlaCasCesTha20,AlvFeh11,IseKroSin19a,IseKroSin18,IseKroSin19,BlaCasMur17,BlaCasEsc22}. However, the cost for computing $U(t)$ is, generally, polynomial in $N$ to our knowledge.
In many cases, e.g., in NMR, the cost of solving \eqref{eq:ODE} can be reduced by compressing the size of the matrix $H(t)$ (e.g., state-space restriction \cite{MenVegPae17,Kup08,KupWagHor07}) or by diminishing the cost of the matrix-vector multiplication (e.g., by exploiting the underlying Kronecker structure \cite{Kup23,AllConLalKup19}). However, the overall cost still increases polynomially in the number of spins \cite{Con20}, and the bottleneck of these approaches remains the solution of an ODE system like the one in \eqref{eq:ODE} \cite{Kup23}. Therefore, deriving a method for $U(t)$ with computation cost that scales linearly with $N$ would mean allowing the simulation of many-body systems unachievable at the moment, i.e., simulating NMR experiments involving more spins than the current state of the art methods achieves.

Recently, a new approach to the problem has been introduced based on the so-called $\star$-product \cite{GisPozInv19,GiPo20,GiPo21}. The idea is to move the problem into a particular algebraic structure of distributions \cite{Ryc23} where it becomes linear. Then the problem is mapped into the usual matrix algebra, truncated, and solved with numerical linear algebra methods \cite{PozVan2023,PozVan2023_leg,PozVan23_mtxeq,Poz23}.
In this paper, we present the first $\star$-approach method competitive with the state-of-the-art. The method is tuned for the case of the generalized Rosen-Zener model \cite{KyoVitSho07}, but it is possible to extend it to other cases, e.g., in the conference proceedings \cite{PozVan23_mtxeq} where a similar approach is used on an example coming from an NMR experiment with Magic Angle Spinning; see also \cite{CiPoRZVB22}.
For the generalized Rosen-Zener model, the cost of the new approach appears to be scaling \emph{linearly} with the size of the problem when computing the operator solution to \eqref{eq:ODE:operator}. The solution is stored in a low-rank format, reducing also the memory cost. The method seems to cheaply scale also as the interval $I$ becomes larger.
Moreover, we present an analogous method for the state vector solution \eqref{eq:ODE}.

After describing the generalized Rosen-Zener model and its use below, the paper recalls the basics of the $\star$-approach (Section \ref{sec:starapproach}). Section \ref{sec:itmethod} introduces the new methods which are tested and compared in Section \ref{sec:exp}. Based on numerical observations, a first cost analysis of the method is presented in Section \ref{sec:cost}.
Section \ref{sec:conclusions} concludes the paper.

\subsection{Generalized Rosen-Zener model}
The Rosen-Zener (RZ) model \cite{RosZen32} is of the highest importance as representative of two-level quantum systems. It corresponds to an exactly soluble system of two coupled time-dependent ordinary differential equations where the driving field is given by a hyperbolic-secant step. The RZ model is of fundamental importance in NMR \cite{SilJosHou85,Hio84} and Magnetic Resonance Imaging (MRI) \cite{ZhaGarPar16}. As a matter of fact, the selective spin inversion in NMR is associated with the exact solution of the so-called Bloch-Riccati equation involving a hyperbolic-secant pulse.   
From the numerical point of view, the RZ model has been central in comparing various approximation methods such as standard perturbation theory, Magnus expansion, Magnus integrators, Euler methods, and explicit Runge-Kutta schemes; see, e.g., \cite{BlaCasOteRos09}. In the framework of quantum-state engineering, Kyoseva et al. \cite{KyoVitSho07,Vit10book} extended the RZ model to multiple degenerate sets of states and higher dimensions of the system of non-autonomous ODE. A few years later, the extended RZ model was used as a test model for solving non-autonomous evolution equations by high-order commutator-free quasi-Magnus exponential integrators \cite{BlaCasTha17} and evaluating the symplectic time average propagators for time-dependent Hamiltonian in the Schrödinger equation \cite{BlaCasMur17}. A posteriori errors for Magnus-type integrators were investigated as well \cite{AuzHofKocQueTha19}. More recently, an efficient algorithm to compute the exponential of skew-Hermitian matrices was propose by Bader et al. \cite{BadBlaCasSey22}.

The dynamics of this extended RZ model is described by the Hamiltonian
\begin{equation}\label{eq:RZ:hamiltonian}
    H(t) := \omega(t) \sigma_3 \otimes I_k + v(t) \sigma_1 \otimes M_k,
\end{equation}
with Pauli matrices
\begin{equation*}
    \sigma_1 = \begin{bmatrix}
        0 & 1 \\
        1 & 0
    \end{bmatrix}, \quad
    \sigma_2 = \begin{bmatrix}
        0 & -i \\
        i & 0
    \end{bmatrix}, \quad
    \sigma_3 = \begin{bmatrix}
        1 & 0 \\
        0 & -1
    \end{bmatrix},
\end{equation*}
$M_k$ is a $k \times k$ tridiagonal symmetric matrix whose main diagonal is composed of all $0$ and the upper and lower diagonals of all $1$, and
where $\omega(t):= w_0 + \varepsilon \cos(\delta t)$, $v(t) := v_0/\cosh(t/T_0)$ are scalar real functions.
Overall, the Hamiltonian has size $N=2k$.

\section{The $\star$-approach to linear ODEs}\label{sec:starapproach}
Consider the system of ODEs
\begin{equation}\label{eq:RZ:ODE}
    \frac{\partial}{\partial t} \psi(t) = -iH(t)\psi(t), \quad \psi(t_0) = \psi_0, \quad t \in I =[t_0,t_f],
\end{equation}
with the Hamiltonian $H(t)$ from \eqref{eq:RZ:hamiltonian}, and the vector $\psi_0\in \mathbb{C}^N$. For the sake of a simpler presentation, assume that $t_0 = -1$ and $t_f = 1$. In fact, any finite interval can be rescaled to $[-1,1]$.

Based on a novel analytical expression \cite{GiLuThJa15,GiPo20,GiPo21} for the solution $\psi(t)$ of \eqref{eq:RZ:ODE}, a new numerical approach for computing $\psi(t)$ was outlined in \cite{PozVan23_mtxeq} by extending the results in \cite{PozVan2023_leg}.
This approach is based on representing the matrix $H(t)\Theta(t-s)$ into a basis of orthonormal Legendre polynomials, where $$\Theta(t-s) = \begin{cases}
    0,\quad \text{if }t<s,\\
    1,\quad \text{if }t\geq s \end{cases},$$ is the \emph{Heaviside theta function}.
The multiplication by $\Theta(t-s)$ is a prerequisite for the use of the novel analytical expression for the solution $\psi(t)$.
We omit this expression here, but we will state the equivalent matrix expression for the Fourier coefficients of $\psi(t)$ in Legendre basis in Equation \eqref{eq:psiApprox} below.

The \emph{coefficient matrix} $F$ of a function $f(t,s)$, is the infinite matrix composed of the Fourier coefficients in the orthonormal Legendre basis $f_{k,\ell} = \int_{-1}^1 \int_{-1}^1 f(t,s) dt \, ds$,
	\begin{equation}\label{eq:coeffMatrix}
		F := \begin{bmatrix}
			f_{k,\ell}
		\end{bmatrix}_{k,\ell=0}^\infty = \begin{bmatrix}
		f_{0,0} & f_{0,1} & f_{0,2} & \dots \\
		f_{1,0} & f_{1,1} & f_{1,2} & \dots \\
		f_{2,0} & f_{2,1} & f_{2,2} & \dots \\
		\vdots & \vdots & \vdots & \ddots
	\end{bmatrix}.
	\end{equation}
The basis of orthonormal Legendre polynomials $\{p_k\}_{k\geq 0}$ is represented by the infinite vector
 \begin{equation}
     \phi(\tau):= \begin{bmatrix}
			p_0(\tau)\\
			p_1(\tau)\\
			\vdots
		\end{bmatrix}.
 \end{equation}
 Thus, a function $f(t,s)$ can be represented in the Legendre basis as
 \begin{equation}\label{eq:infDoubleSeries}
     f(t,s) = \sum_{k,\ell} f_{k,\ell} p_k(t)p_{\ell}(s) = \phi(t)^\top F \phi(s).
 \end{equation}
 In case the function $f(t,s)$ is of the form $f(t,s) = \Tilde{f}(t) \Theta(t-s)$, with $\Tilde{f}(t)$ a continuous function, then the equality holds except for $t=s$.

 In practice, we will work with a finite truncation\footnote{The nontrivial analysis of the truncation error can be found in \cite{PozVan2023_leg} for the case $N=1$. Until now, all the numerical experiments show that the same results hold also for systems of linear ODEs, i.e., $N>1$. A rigorous extension of this analysis to systems of ODEs is out of the scope of this paper and will be developed in future work.} of the double series \eqref{eq:infDoubleSeries}, i.e., 
 \begin{equation}\label{eq:finDoubleSeries}
     f(t,s) \approx \sum_{k,\ell=0}^{M-1} f_{k,\ell} p_k(t)p_{\ell}(s) =: \phi_M(t)^\top F_M \phi_M(s),
 \end{equation}
 where $\phi_M(\tau) := \begin{bmatrix}
     p_0(\tau) & p_1(\tau) &\dots & p_{M-1}(\tau)
 \end{bmatrix}^\top$ and $F_M$ is the $M\times M$ leading principal submatrix of $F$ \eqref{eq:coeffMatrix}. 
Let $T_M$ denote the $M\times M$ leading principal submatrix of the coefficient matrix for $\Theta(t-s)$, then we know from \cite{PozVan23_mtxeq} that 
\begin{equation}\label{eq:psiApprox}
    \psi(t) \approx (I_N\otimes \phi(t)^\top T_M) x,  \quad t \in [-1, 1],
\end{equation}
where $x\in\mathbb{C}^{MN}$ is the solution to the linear system of equations
\begin{equation}\label{eq:sysEq}
    (I_{MN}+i\mathcal{H}_{MN}) x = \psi_0 \otimes \phi_M(-1),
\end{equation}
with $\mathcal{H}_{MN}$ the block matrix obtained by representing each element in the matrix $H(t)\Theta(t-s)$ by its (truncated) coefficient matrix and $\otimes$ the Kronecker product.
That is, let $W_M\in\mathbb{C}^{M\times M}$ denote the coefficient matrix of $\omega(t)\Theta(t-s)$ and $V_M\in\mathbb{C}^{M\times M}$ the coefficient matrix of $v(t)\Theta(t-s)$, then we have
\begin{align*}
    \mathcal{H}_{MN} &= (\sigma_3 \otimes I_k) \otimes \Omega_M + (\sigma_1 \otimes M_k) \otimes V_M\\
    &= \left[
				\begin{array}{cccc|cccc}
					\Omega_M & & & & 0 & V_M\\
					& \Omega_M & & & V_M & 0 & \ddots \\
					& & \ddots  & & & \ddots & \ddots & V_M\\
					& &  & \Omega_M  & & & V_M & 0\\
					\cline{1-8}
					0& V_M & & & -\Omega_M \\
					V_M& 0 & \ddots & & & -\Omega_M\\
					& \ddots& \ddots & V_M & & & \ddots \\
					& & V_M& 0 & & & & -\Omega_M
				\end{array}\right].
\end{align*}
The matrix structure can be exploited to obtain efficient solvers for the system of equations \eqref{eq:sysEq}.

\section{A new iterative method for the Rosen-Zener ODE}\label{sec:itmethod}
To fully exploit the Kronecker structure of the matrix $\mathcal{H}_{MN}$, we will reformulate \eqref{eq:psiApprox} as
\begin{equation}\label{eq:approx:psi}
    \psi(t) \approx vec\left(\phi_M(t)^T T_M X\right), \quad t \in [-1, 1],
\end{equation}
where $vec$ is the vectorization transformation, and $X \in \mathbb{C}^{M \times N}$ is the solution of the matrix equation reformulation of \eqref{eq:sysEq}, i.e., 
\begin{equation}\label{eq:RZ:mtxeq}
    X + i\,\Omega_M X (\sigma_3 \otimes I_k) + i\,V_M X (\sigma_1 \otimes M_k) = \phi_M(-1) \psi_0^\top;
\end{equation}
for an introduction to matrix equations and their numerical solution see, e.g., \cite{Sim16}.
In other words, the solution of the linear system of equations \eqref{eq:sysEq} is given by $x = vec(X)$.
One advantage of the matrix equation formulation is that it uses smaller matrices, and thus, is more memory efficient.

\bigskip

In order to solve \eqref{eq:RZ:mtxeq} we make use of the (implicit) iterates:
\begin{equation*}
    X_{n+1} + i\,\Omega_M X_{n+1} (\sigma_3 \otimes I_k) = -i\,V_M X_{n} (\sigma_1 \otimes M_k) + \phi_M(-1) \psi_0^T.
\end{equation*}
Thanks to the simple diagonal structure of the matrix $\sigma_3 \otimes I_k$, these iterates can be rewritten as the following \emph{stationary iterative method} (fixed point method)
\begin{align}
    \label{eq:iterates1} X_{n+1/2} &=  -i\,V_M X_n (\sigma_1 \otimes M_k) + \phi_M(-1) \psi_0^T; \\
    \label{eq:iterates2} X_{n+1}   &=  G_1 X_{n+1/2} D_1 + G_2 X_{n+1/2} D_2,
\end{align}
with $X_0 = \phi_M(-1) \psi_0^T$ and
\begin{align*}
    G_1 &= (I_M +i\, \Omega_M)^{-1}, \quad D_1 = \begin{bmatrix}
        1 & 0 \\
        0 & 0
    \end{bmatrix} \otimes I_k, \\
    G_2 &= (I_M -i\, \Omega_M)^{-1}, \quad D_2 = \begin{bmatrix}
        0 & 0 \\
        0 & 1
    \end{bmatrix} \otimes I_k.
\end{align*}
The iterates \eqref{eq:iterates1}--\eqref{eq:iterates2} can also be vectorized, being transformed into the following ones
\begin{align}\label{eq:it}
    x_{n+1} &= -G(i\, (\sigma_1\otimes M_k)\otimes V_M) x_n + G(\psi_0 \otimes \phi_M(-1)), 
    \end{align}
    with $G = (I_{NM}+i\,(\sigma_3 \otimes I_k)\otimes\Omega_M)^{-1}$.
From classical results on stationary iterative methods (e.g., \cite{saad03}), the method converges if the spectral radius of the iteration matrix is smaller than $1$, i.e., $\rho(G(i\, (\sigma_1\otimes M_k)\otimes V_M))<1$. Moreover, the smaller $\rho(G(i\, (\sigma_1\otimes M_k)\otimes V_M))$ is, the faster the (linear in log-scale) convergence.

Naturally, the latter vectorized expressions have only a theoretical purpose. Indeed, the iterates \eqref{eq:iterates1}-\eqref{eq:iterates2} are computationally less expensive. Moreover, the matrix expression $X$ of the solution allows us to reduce the computation cost further, as we explain in the next section.
    
\subsection{Low-rank approach}
As noticed in \cite{PozVan23_mtxeq}, the solution $X$ of the matrix equation \eqref{eq:RZ:mtxeq} is often characterized by a low numerical rank, that is, given the singular value decomposition (SVD) $X = USV^H$ many of the singular values of $X$, $s_1 \geq s_2\geq \dots \geq s_{\min(M,N)}$, are close or equal to zero.
More specifically, the $s_j$ tend to decay exponentially to zero as $j$ increases. This \emph{decay phenomenon} is well-known in literature, and it is often associated with matrix equations with a low-rank right-hand side, as $\phi_M(-1)\psi_0^T$ is in our case; see, e.g., \cite{Sim16}.

This means that we can try to approximate the solution $X$ by SVD truncation, i.e., setting to zero all the $s_r, \dots, s_{\min(M,N)}$ smaller than a given tolerance. As a consequence, we get the approximation $X \approx L R^T$, with $L = U(1:M,1:r)S(1:r,1:r) \in \mathbb{C}^{M\times r}, R = V(1:N,1:r) \in \mathbb{C}^{N \times r}$, with $r << \min(M,N)$ (note that we are using Matlab notation, where $A(i:j,k:\ell)$ denotes the submatrix formed by the $i$th through $j$th row and $k$th through $\ell$th column of $A$).  Moreover, starting from the rank-$1$ matrix $X_0 = \phi_M(-1)\psi_0^T$, we can therefore build a method that tries to produce a low-rank approximation for each iterate, i.e., $X_n \approx L_n R_n^T$. As a consequence, given general matrices $A, B$, we can approximate the matrix product $A X_n B^T$ by the cheaper product $A L_n (B R_n)^T$. This low-rank approach is nowadays standard in many matrix equation solvers; see, e.g., \cite{Sim16}.
Combining iterates \eqref{eq:iterates1}--\eqref{eq:iterates2} with the described low-rank approach, we obtain the new Algorithm \ref{alg:vec}.

\begin{algorithm}
\caption{Iterative low-rank method for Equation \eqref{eq:RZ:mtxeq}}\label{alg:vec}
\begin{algorithmic}
\Require Error tolerance $\textsc{tol} > 0$ and svd truncation tolerance $\textsc{trunc} > 0$.
\State $L = \phi_M(-1); \quad R = \psi_0;$
\State $g =  [\, G_1 L,\, G_2 L \,]; \quad d = [\,D_1 R,\, D_2 R\,];$
\While{$err\geq \textsc{tol}$}
    \State $L = -i\, V_M L; \quad\quad\quad\quad R = (\sigma_1 \otimes M_k) R;$  \Comment{Iteration \eqref{eq:iterates1}}
    \State $L =  [\,G_1 L,\, G_2 L,\, g\,]; \quad R = [\,D_1 R,\, D_2R,\, d\,];$ \Comment{Iteration \eqref{eq:iterates2}}
    \State $L = Q_L R_L;$   \Comment{Economy-size \;QR\; decomposition}
    \State $R_L = U S V^H;$ \Comment{Economy-size SVD decomposition}
    \State  $r = \min\{j : S(j,j) < \textsc{trunc} \}$;
    \State $L = Q_L \, U(:,1:r) \, S(1:r,1:r);$   \Comment{Truncation}
    \State $R = R  \, \text{ conjugate}(V(:,1:r));$ \Comment{Truncation}
    \State $b = L (R^T \text{ conjugate}(\psi_0));$  \Comment{$b = (\psi_0^H\otimes I_M) vec(L_{n} R_{n}^T)$}
    \State $err = \| b - b_{old} \|_2 $     \Comment{Cheap error estimate}
    \State $b_{old} = b$;
\EndWhile
\end{algorithmic}
\end{algorithm}

The SVD in Algorithm \ref{alg:vec} is applied only to the left-hand factor $L$, since, in our numerical experiments, we noticed that the low-rank property of the solution seems to be associated with the left-hand side of the equation. Moreover, if the size $N$ of the ODE is large, the singular values decomposition of the right-hand side $R$ becomes too expensive.

Note that the stopping criterion is based on the idea of computing the error estimate $\|b - b_{old} \|_2 = \| (\psi_0^H\otimes I_M) vec(L_{n}R_{n}^T) - (\psi_0^H\otimes I_M) vec(L_{n-1}R_{n-1}^T) \|_2$ for the quantity $b =(\psi_0^H\otimes I_M) vec(X_n)$. This is a \emph{cheap} estimate of the error since it avoids computing the matrix $X_n$ from the factors $L_n, R_n$. If needed, $\psi_0$, in the expression for $b$, can be replaced by other nonzero vectors.
 
\subsection{Computing the operator solution}
Consider now the problem of computing the operator solution, i.e., the $N \times N$ matrix-valued function $U(t)$ solving
\begin{equation}\label{eq:RZ:ODE:operator}
    \frac{\partial}{\partial t} U(t) = -iH(t)U(t), \quad U(t_0) = I_N, \quad t \in I =[t_0,t_f].
\end{equation}
To compute $U(t)$ we can solve Equation~\eqref{eq:RZ:mtxeq} with $\psi_0 = e_j$, for $j=1,\dots N$, denoting the solution as $X^{(j)}$. Then, using approximation \eqref{eq:approx:psi} we get
\begin{equation*}
    U(t)e_j \approx vec\left(\phi_M(t)^T T_M X^{(j)}\right).
\end{equation*}
In order to approximate $X^{(j)}$, we can use Algorithm \ref{alg:vec}. Moreover, we can combine the $N$ runs of the algorithm into one algorithm, allowing the low-rank approximation of $\mathcal{X} = [X^{(1)}, \dots, X^{(N)}]$. In this way, we can perform just one singular value truncation for the $M \times N^2$ matrix $\mathcal{X}$ per iteration. This idea results in Algorithm \ref{alg:mtx} where we marked in blue the main changes with respect to Algorithm \ref{alg:vec}.
Note that the columns of the approximated $\mathcal{X}$ in Algorithm \ref{alg:mtx} are ordered differently, that is, the factors $L$ and $R$  are so that 
\begin{equation*}
    X^{(j)} \approx L \, \big(R(1:N,j:N:j+(r-1)N)\big)^T, \quad j=1,\dots, N. 
\end{equation*}

\begin{algorithm}
	\caption{Iterative low-rank method for the operator solution}\label{alg:mtx}
	\begin{algorithmic}
		\Require Error tolerance $\textsc{tol} > 0$ and svd truncation tolerance $\textsc{trunc}> 0$.
		\State $L = \phi_M(-1); \quad {\color{blue}R = I_N};$
		\State $g =  [\, G_1 L,\, G_2 L \,]; \quad d = [\,D_1,\, D_2\,];$
		\While{$err\geq \textsc{tol}$}
		\State $L = -i\, V_M L; \quad\quad\quad\quad R = (\sigma_1 \otimes M_k) R;$  \Comment{Iteration \eqref{eq:iterates1}}
		\State $L =  [\,G_1 L,\, G_2 L,\, g\,]; \quad R = [\,D_1 R,\, D_2R,\, d\,];$ \Comment{Iteration \eqref{eq:iterates2}}
		\State $L = Q_L R_L;$   \Comment{Economy-size \;QR\; decomposition}
		\State $R_L = U S V^H;$ \Comment{Economy-size SVD decomposition}
		\State  $r = \min\{j : s_j < \textsc{trunc} \}$, with $\text{diag}(s_j) = S$;
		\State $L = Q_L \, U(:,1:r) \, S(1:r,1:r);$   \Comment{Truncation}
		\State {\color{blue}$K = \text{conjugate}(V(:,1:r)) \otimes I_N;$}
		\State {\color{blue}$R = R K;$} \Comment{Truncation}
		\State {\color{blue}$b = L (R(1,1:N:1+(r-1)N))^T;$}  \Comment{$b = \mathcal{X}(1:M,1)$}
		\State $err = \| b - b_{old} \|_2 $     \Comment{Cheap error estimate}
		\State $b_{old} = b$;
		\EndWhile
	\end{algorithmic}
\end{algorithm}

Note that the stopping criterion is based on an error estimate for the approximant of the first column of $\mathcal{X}$. This criterion is cheap but clearly naive; however, it is pretty effective in all the presented numerical experiments. For this reason, testing and developing more mathematically founded criteria is out of the scope of this paper.

\section{Numerical experiments and comparisons}\label{sec:exp}
Following \cite{BlaCasMur17}, in the Rosen-Zener model \eqref{eq:RZ:hamiltonian} we set $w_0 = 5$ and $v_0 = 1/2$ and we consider the cases:
\begin{enumerate}[(a)]
    \item $\varepsilon = 0$, $T_0 = 10$;
    \item $\varepsilon = 1/10$, $\delta = 1/10$, $T_0 = 5$;
    \item $\varepsilon = 1/2$, $\delta = 1$, $T_0 = 5$;
    \item $\varepsilon = 2$, $\delta = 5$, $T_0 = 1$.
\end{enumerate}
In this section, we compare the proposed algorithms with standard Runge-Kutta (RK) methods and with methods\footnote{The methods have been implemented in MatLab by the authors of \cite{BlaCasMur17} and released at \url{https://www.gicas.uji.es/Research/TD-propagators.html}. In order to make a fair comparison, we implemented our own version of them that exploits MatLab sparse matrices.} that have been tested  in \cite{BlaCasMur17} on the same cases we are considering. They are:
\begin{itemize}
    \item SM$_8^{[4]}$, SM$_{11}^{[6]}$, SM$_{11}^{[8]}$: Respectively the $8$-stage $4$th-order, $11$-stage $6$th-order, and $11$-stage $8$th-order symplecting splitting methods introduced in \cite{BlaCasMur17};
    \item S$_6^{[4]}$, S$_{10}^{[6]}$: Respectively, the $6$-stage $4$th order and the $10$-stage $6$th-order methods from \cite{BlaMoa02} adapted to this case as described in \cite{BlaCasMur17};
    \item RK$_7^{[6]}$: the $7$-stage $6$th-order explicit Runge-Kutta method (with Lobatto quadrature rule) from \cite{BlaCasMur17}.
\end{itemize}
After a first experiment meant to show the accuracy of Algorithm \ref{alg:vec}, we focus on the operator solution \eqref{eq:RZ:ODE:operator}. In this latter case, we refer by the \emph{$\star$-approach} (denoted by ``star'') to the method obtained by computing the coefficient matrices $\Omega_M, V_M$ and $T_M$ (discretization) and solving the related matrix equations by Algorithm \ref{alg:mtx}.
The experiments have been implemented in MatLab and have been obtained by running MATLAB R2022a on a laptop with Intel I7 CPU.
In Experiment 2, 3 and 4 we used the $8$th-order Magnus integrator from \cite{BlaCasRos00} as the reference solution for the error computation, which is also the reference solution in \cite{BlaCasMur17}.

\begin{figure}[htbp]
    \centering
    \begin{subfigure}{\textwidth}
    \includegraphics[width=0.43\textwidth]{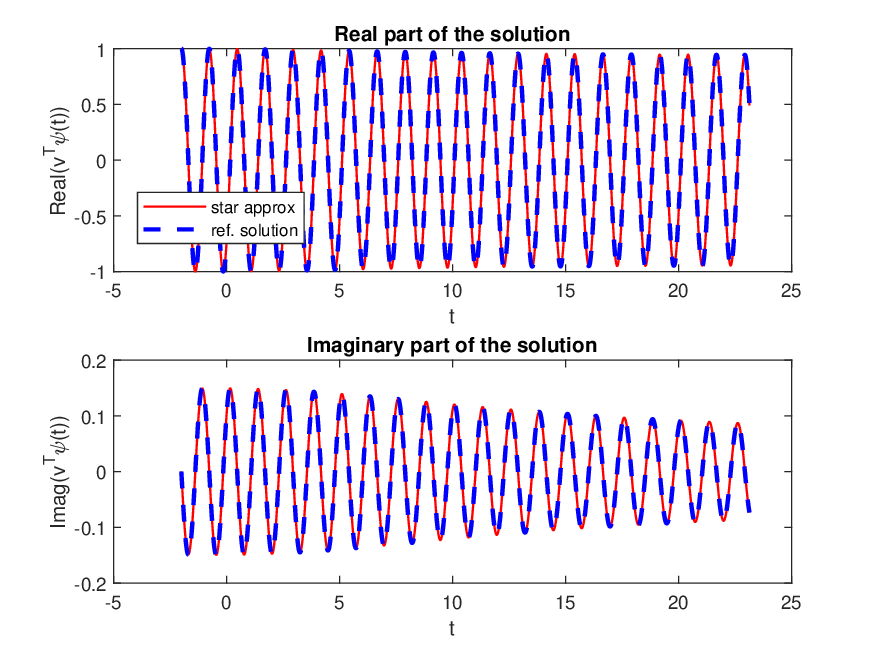}
    \includegraphics[width=0.43\textwidth]{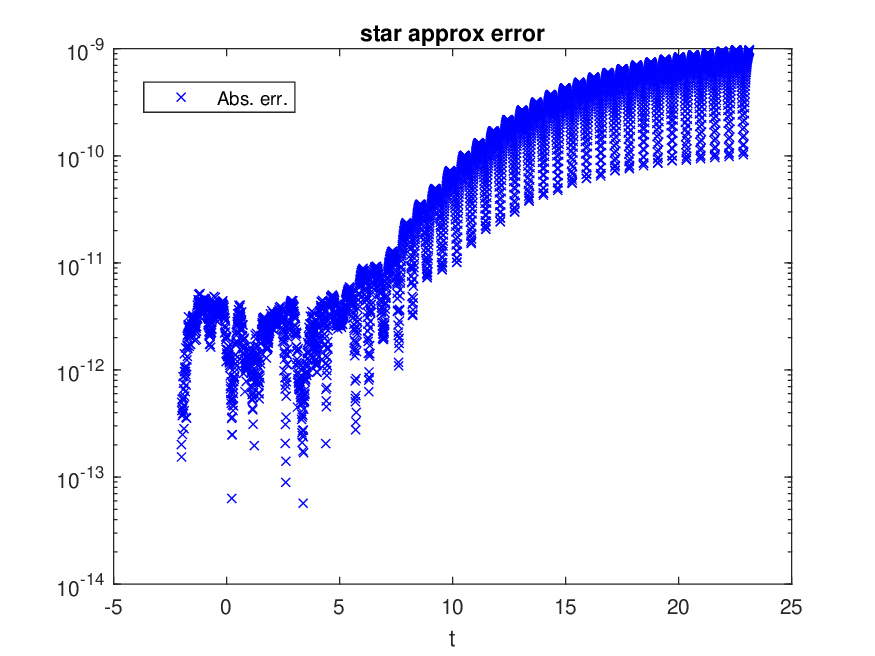} 
    \caption{\;\;\;\;\;\;\;\;\;\;\;\;\;\;\;\;}
     \end{subfigure}
\begin{subfigure}{\textwidth}
    \includegraphics[width=0.43\textwidth]{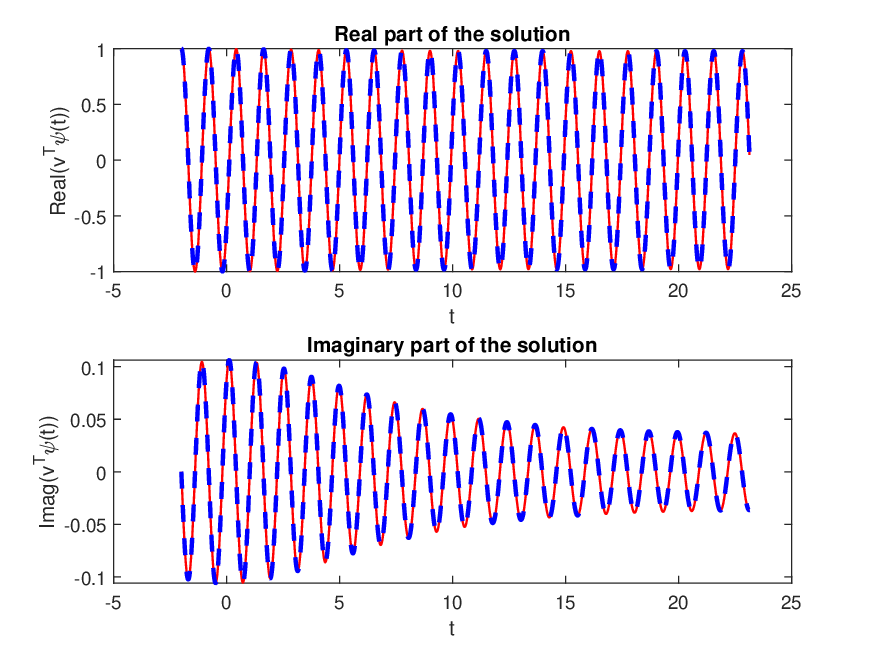}
    \includegraphics[width=0.43\textwidth]{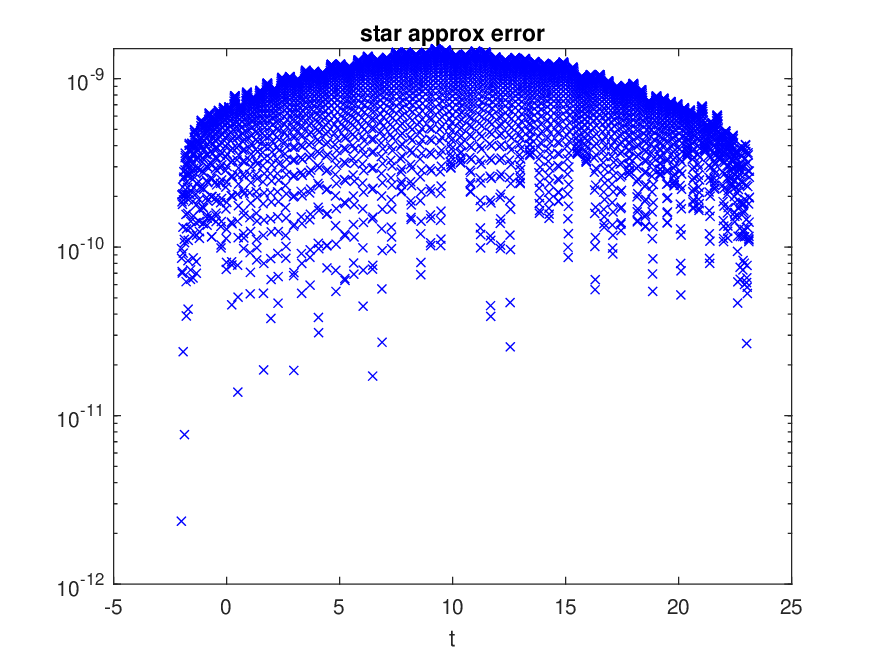} 
    \caption{\;\;\;\;\;\;\;\;\;\;\;\;\;\;\;\;}
     \end{subfigure}
     \begin{subfigure}{\textwidth}
    \includegraphics[width=0.43\textwidth]{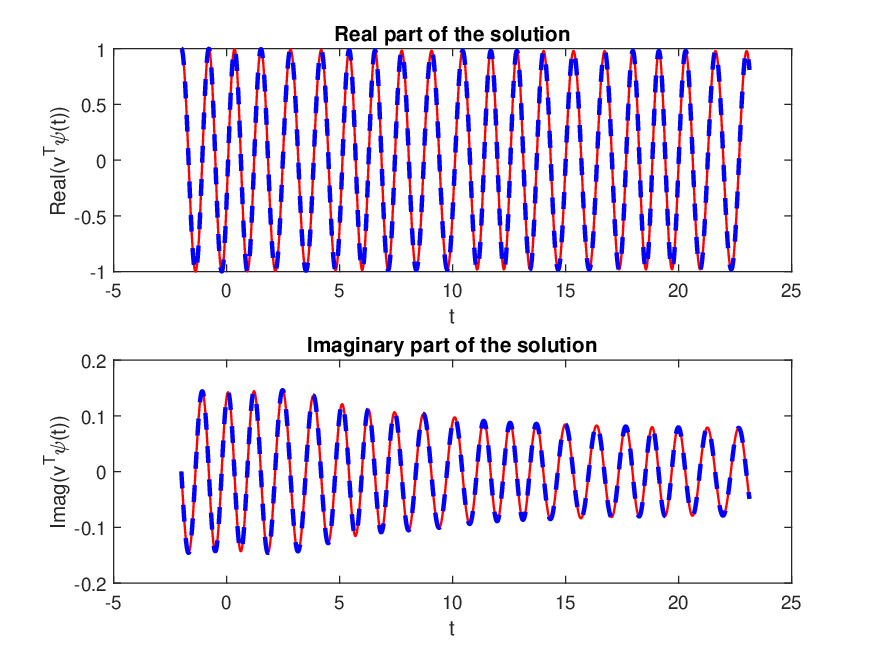}
    \includegraphics[width=0.43\textwidth]{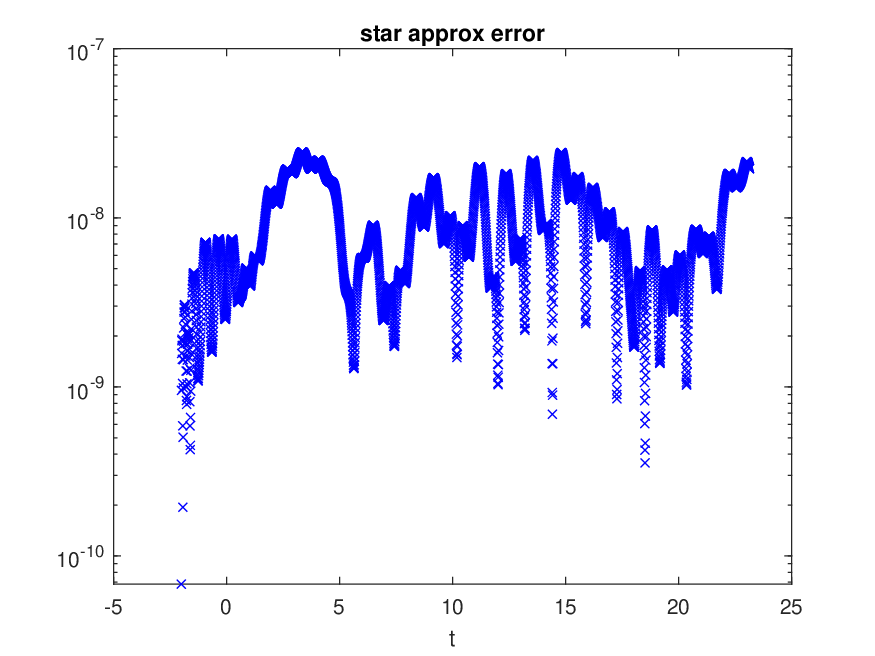} 
    \caption{\;\;\;\;\;\;\;\;\;\;\;\;\;\;\;\;}
     \end{subfigure}
     \begin{subfigure}{\textwidth}
    \includegraphics[width=0.43\textwidth]{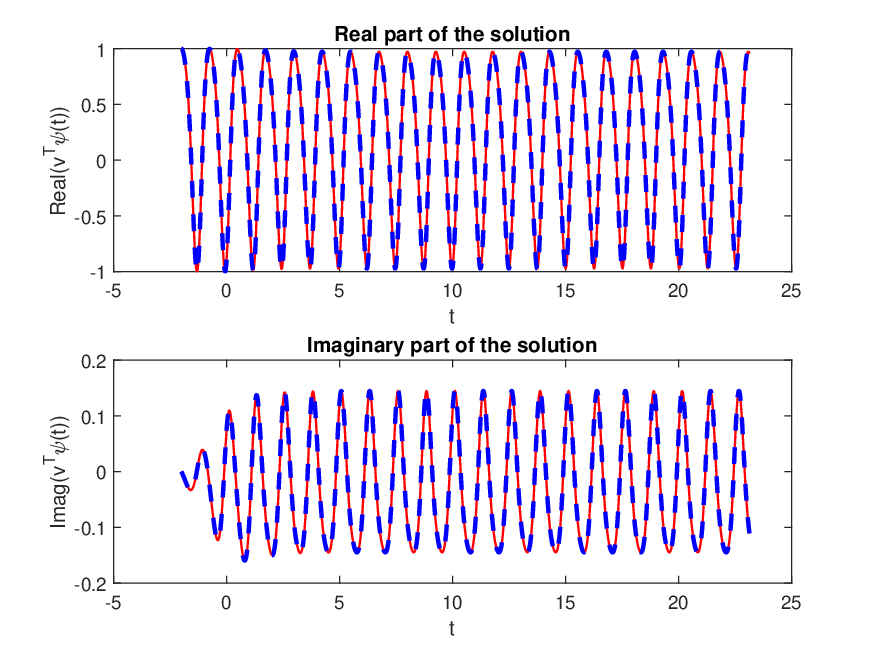}
    \includegraphics[width=0.43\textwidth]{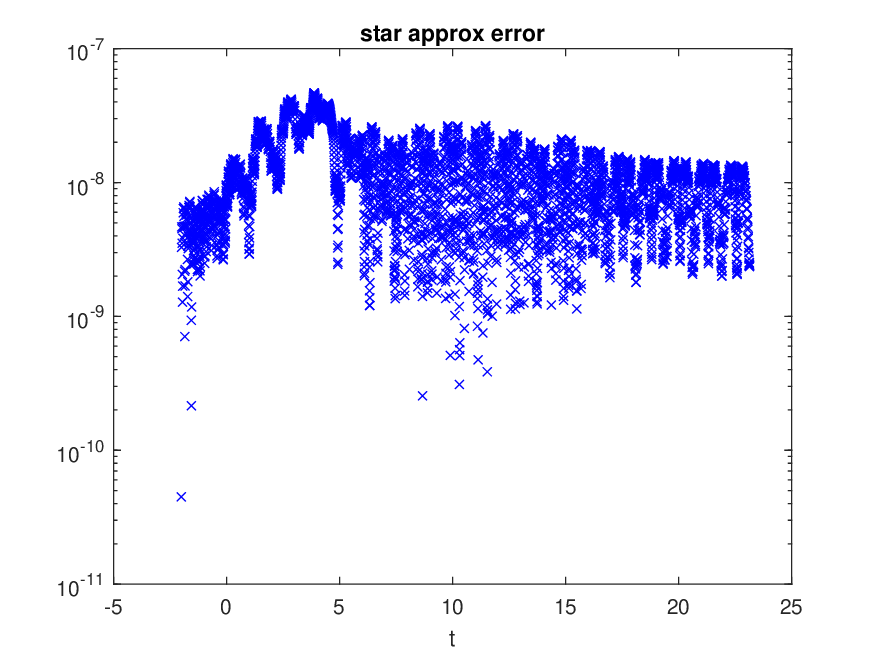} 
    \caption{\;\;\;\;\;\;\;\;\;\;\;\;\;\;\;\;}
     \end{subfigure}
    \caption{Experiment 1: State vector ODE \eqref{eq:ODE} for cases (a)-(d) with a normalized random starting state $\psi_0$. Left: Approximation obtained by Algorithm \ref{alg:vec} (red) and reference solution (blue dashed). Right. Absolute error of the approximation to the reference solution.}
    \label{fig:psi}
\end{figure}

\paragraph{Experiment 1}
In the first experiment, we test the accuracy of Algorithm \ref{alg:vec}, proving that it is able to solve the ODEs of interest.
For each of the four cases, we approximate the value $\psi_0^H \psi(t)$  with initial time $t_0 = -2$, final time $t_f = t_0 + 8\pi$ and starting state $\psi_0$ as a normalized random vector.
Using Algorithm \ref{alg:vec} and \eqref{eq:approx:psi}, we obtain the approximation
\begin{equation*}
    \psi_0^H \psi(t) \approx \beta(t) := \phi_M(t)^T T_M L (\psi_0^H R)^T. 
\end{equation*}
The tolerance for the stopping criterion of Algorithm \ref{alg:vec} is set to $\textsc{tol}=1\cdot10^{-7}$ and the SVD truncation tolerance to $\textsc{trunc}=1\cdot 10^{-6}$. 
As a reference for the error estimates, we consider the Runge-Kutta $(4,5)$ formula implemented in the Matlab function \verb|ode45| with absolute and relative tolerances set to $3\cdot 10^{-12}$. Then we compute $\beta(t)$ on the time points given by \verb|ode45|, so obtaining the error plots in Figure \ref{fig:psi}. Note that the error in the plots is always below the tolerance $\textsc{tol}$. In Table \ref{table:params}, we report the settings, the number of iterations, and the maximal number of singular values $r$ kept by Algorithm \ref{alg:vec}. Observe that the number of iterations is limited, and the maximal number of singular values is always much smaller than $M$.

\begin{table}[!ht]
\centering
\begin{tabular}{ |c||c|c|c|c|  }
 \hline
 Case & $M$ & max err & n it & max n sing val\\
 \hline
 (a) & $130$ & $9.7788 \cdot10^{-10}$ & $26$ & $33$ \\
 (b) & $140$ & $1.5059\cdot10^{-9}$ & $22$ & $33$ \\
 (c) & $250$ & $2.4463\cdot10^{-8}$ & $22$ & $41$ \\
 (d) & $550$ & $4.6723\cdot10^{-8}$ & $15$ & $53$ \\
 \hline
\end{tabular}
\caption{Experiment 1: State vector ODE \eqref{eq:ODE} for cases (a)-(d) with a normalized random starting state $\psi_0$. For Algorithm \ref{alg:vec}, from left to right: truncation parameter $M$, maximal absolute error on the whole interval $ [t_0,t_f] $, number of iterations, and maximal number of singular values $r$.}
\label{table:params}
\end{table}

\begin{figure}[htbp]
    \centering
    \begin{subfigure}{\textwidth}
    \includegraphics[width=0.5\textwidth]{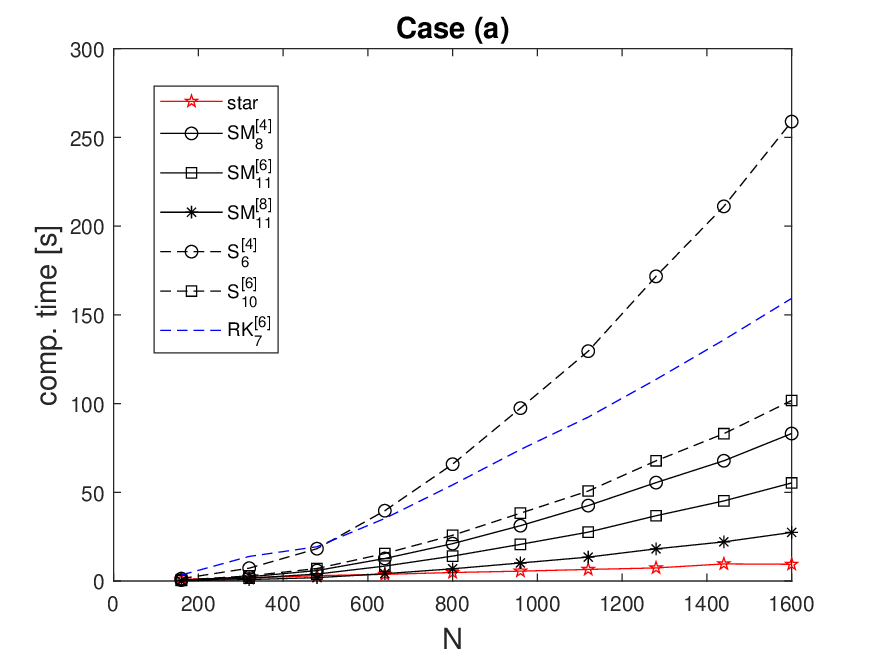}
    \includegraphics[width=0.5\textwidth]{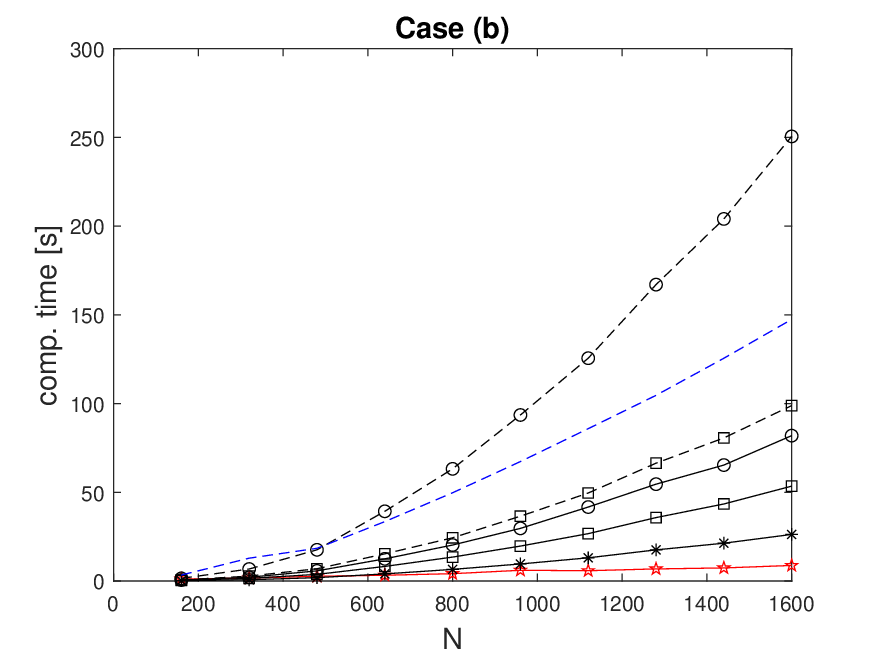} 
     \end{subfigure}
\begin{subfigure}{\textwidth}
    \includegraphics[width=0.5\textwidth]{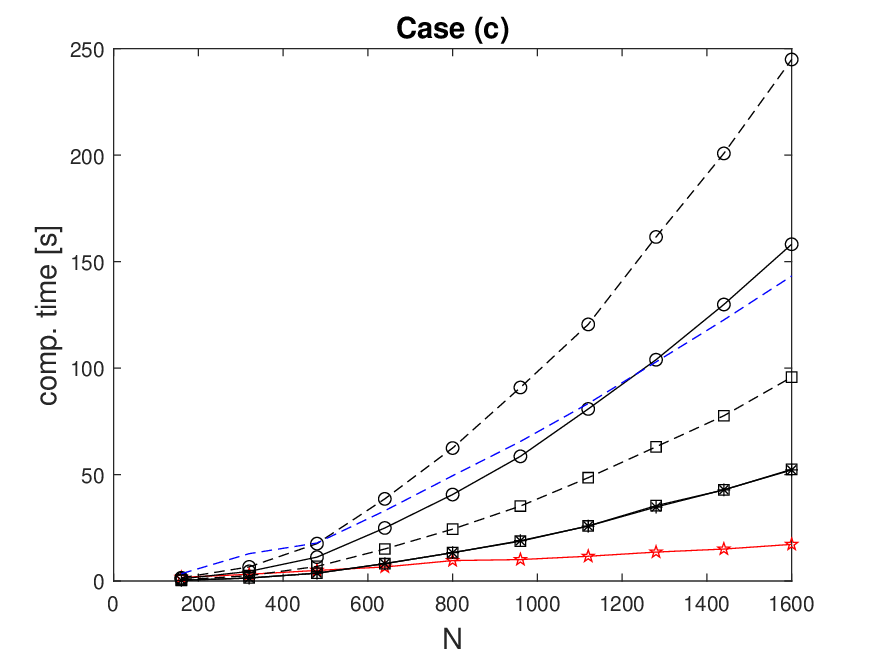}
    \includegraphics[width=0.5\textwidth]{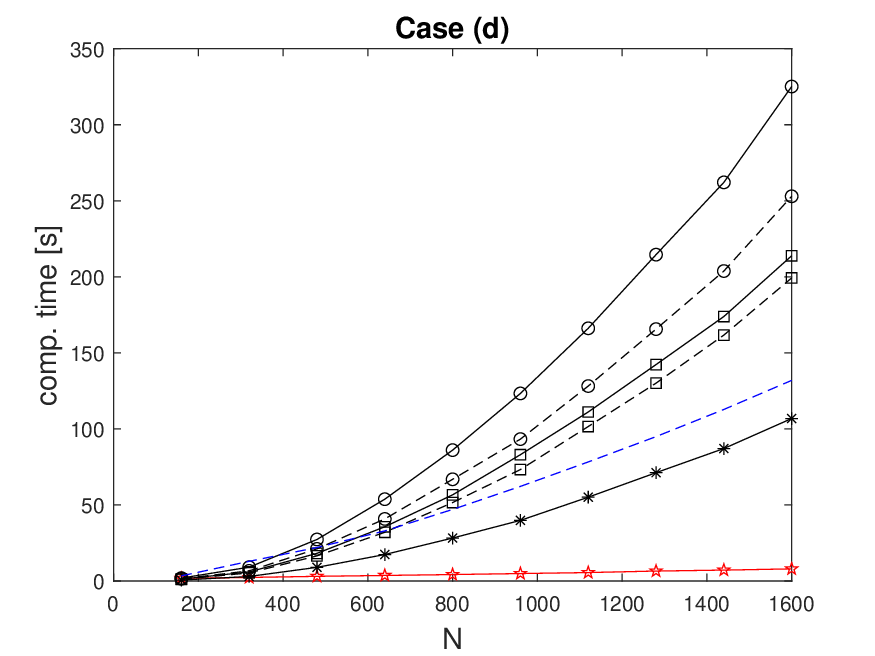} 
     \end{subfigure}
    \caption{Experiment 2. computation time comparison on the operator solution of cases (a)--(d) for increasing system size $N$. All the methods reach an absolute error (Euclidean matrix) norm smaller than $1\cdot 10^{-6}$ at the finale time $t_f$.}
    \label{fig:ctime}
\end{figure}

\paragraph{Experiment 2}
In this experiment, we test the computation time of the $\star$-approach (star) for the operator solution of cases (a)--(d) as the size of the system $N$ increases. 
The computation time is obtained by summing the time of Algorithm \ref{alg:mtx} (solving the matrix equation) and of the discretization (computing the coefficient matrices).
We compare the algorithm with the methods described above. In Figure~\ref{fig:ctime}, we observe that the computation time of the $\star$-approach appears to be linear, or at least sub-quadratic, as $N$ increases, while all the other methods' computation time grows quadratically. This is particularly evident in case (d). In general, while the $\star$-approach is not competitive for small $N$, it becomes extremely competitive as $N$ increases. All the methods' parameters have been set to be the optimal in order to reach an absolute error matrix (Euclidean) norm at the final time $t_f$ smaller than $1\cdot 10^{-6}$. The tolerance for the stopping criterion of Algorithm \ref{alg:vec} is set to $\textsc{tol}=1\cdot10^{-7}$ and the SVD truncation tolerance to $\textsc{trunc}=1\cdot 10^{-6}$. The other parameters can be found in Tables \ref{tab:ctime:a} and \ref{tab:ctime:bcd} (Appendix), where $2^\kappa$ is the number of subintervals in which the domain has been split. Table~\ref{tab:ctime:a} reports also the absolute error norm at $t_f$ for each method and for each considered $N$ in case (a). In order to avoid repetitions, Table~~\ref{tab:ctime:bcd} reports the norm of the error for each method only for $N=1600$. The results for the other values of $N$ are analogous.

\begin{figure}[htbp]
    \centering
    \begin{subfigure}{\textwidth}
    \includegraphics[width=0.5\textwidth]{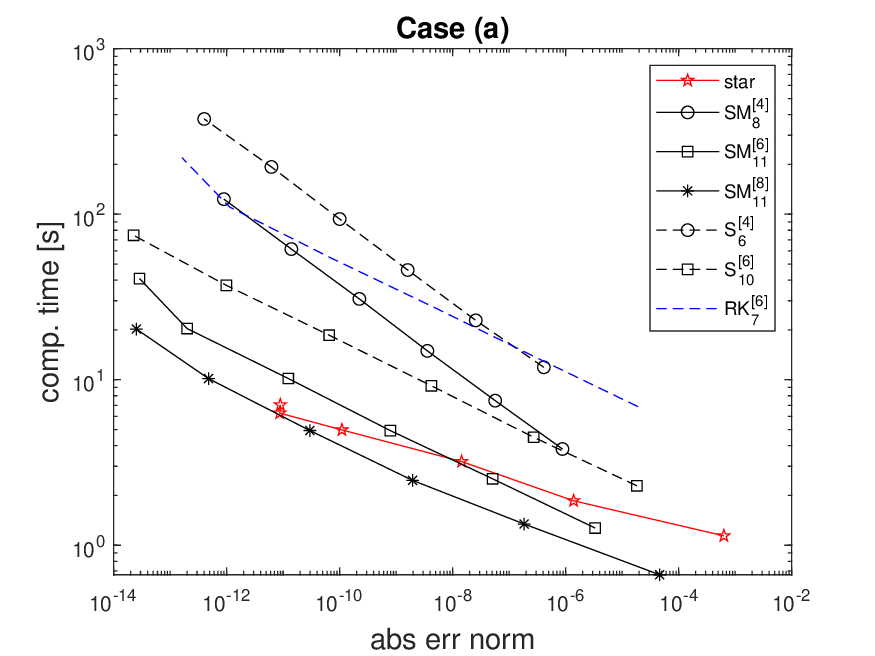}
    \includegraphics[width=0.5\textwidth]{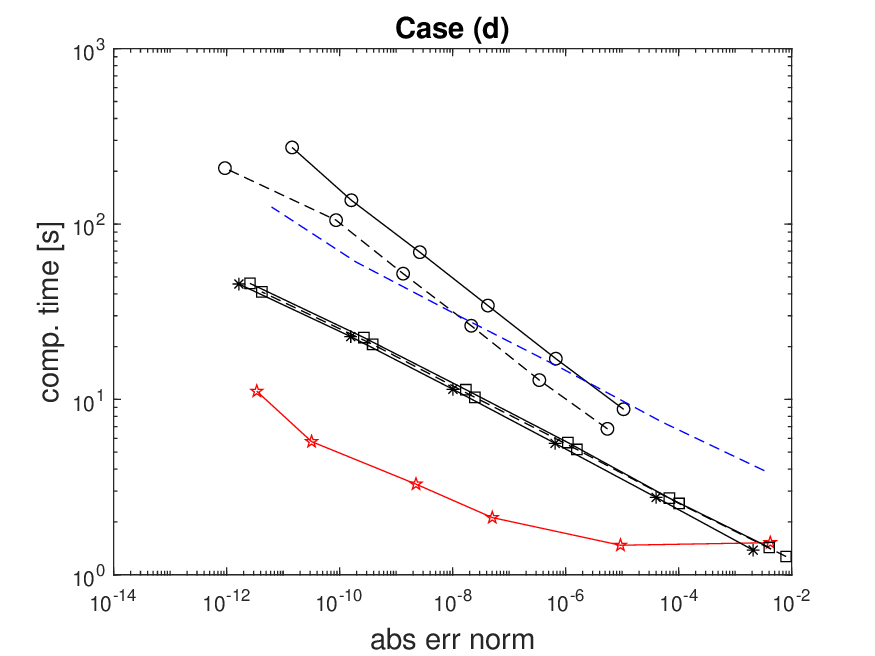} 
     \end{subfigure}
    \caption{Experiment 3. Computation time comparison on the operator solution of cases (a) and (d) for increasing absolute error norm accuracy. The system size is $N=400$.}
    \label{fig:accuracy}
\end{figure}

\paragraph{Experiment 3}
We compare the computation time needed by each method to reach different levels of accuracy for the operator solution described in Experiment 2. Figure~\ref{fig:accuracy} reports the results for cases (a) and (d), for $N=400$.  As we can see, the $\star$-approach (star) compares well with the other methods. It is important to remark that Algorithm \ref{alg:mtx} cannot reach higher accuracy than the one in the plots as the iterations stagnate. A possible explanation is that the linear system \eqref{eq:sysEq} is too ill-conditioned to obtain higher accuracy.

\begin{figure}[htbp]
    \centering
    \begin{subfigure}{\textwidth}
    \includegraphics[width=0.5\textwidth]{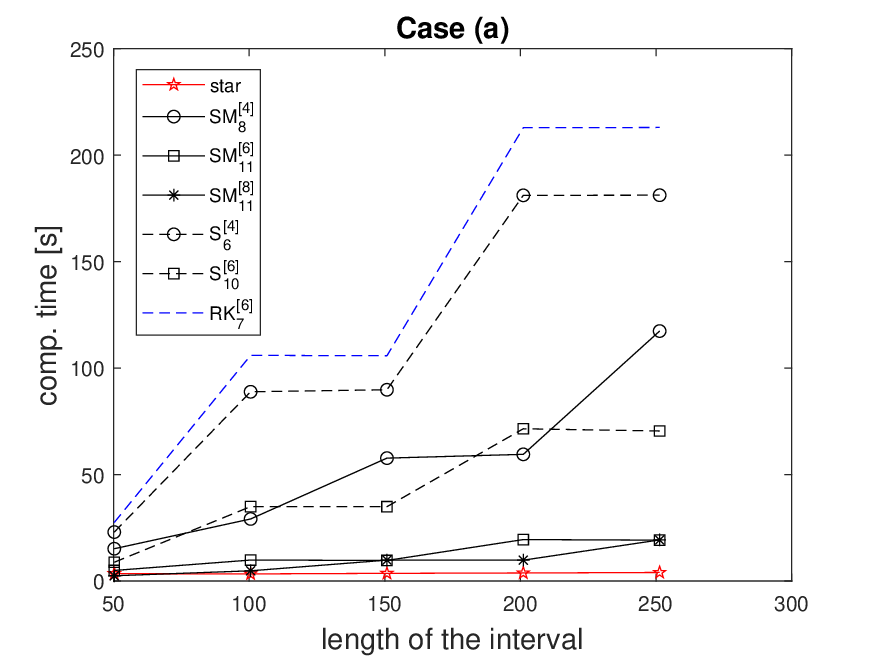}
    \includegraphics[width=0.5\textwidth]{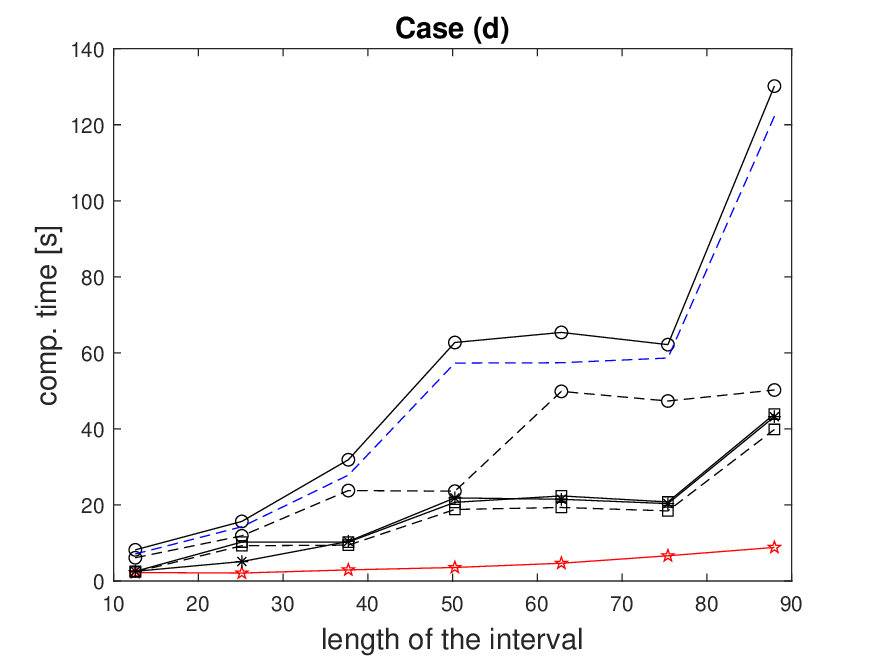} 
     \end{subfigure}
    \caption{Experiment 4. Computation time comparison on the operator solution of cases (a) and (d) for increasing length of the interval $I$. The system
size is $N = 400$.}
    \label{fig:interval}
\end{figure}

\paragraph{Experiment 4}
In this last experiment, we consider the operator cases (a) and (d) with $N=400$ and the length of the interval $I = [t_0,t_f]$ is increased by increasing $t_f$. Figure~\ref{fig:interval} shows the results. Again, the $\star$-approach outperforms the other methods. As before, the computation time of the $\star$-approach is the sum of the time needed by Algorithm \ref{alg:mtx} and the time needed for computing the coefficient matrices (discretization). In these experiments, the computation time for Algorithm \ref{alg:mtx} scales linearly with the length of $I$, while the discretization time is quadratic. In the case (d), for the longest intervals, the discretization time starts dominating the computation time. Tables \ref{tab:interval:a} and \ref{tab:interval:d} in the Appendix report the parameter settings and the absolute error norm at $t_f$ for each method and each interval. In Algorithm \ref{alg:vec}, we set $\textsc{tol}=1\cdot10^{-7}$ and $\textsc{trunc}=1\cdot 10^{-6}$.

\section{A first cost explanation}\label{sec:cost}
Proving that the computation time of the $\star$-approach in Figure~\ref{fig:ctime} scales linearly with $N$ requires proving the following Properties:
\begin{enumerate}
    \item[P1.] The number of iterations of Algorithm \ref{alg:mtx} does not asymptotically increases with $N$;
    \item[P2.] The truncation parameter $M$ does not depend on $N$; 
    \item[P3.] The number of singular values larger than the given tolerance is smaller than $M$ at every iteration for a fixed large enough $M$; 
    \item[P4.] The number of nonzero elements of $R$ at each iteration scales linearly with $N$.
\end{enumerate}
Indeed, first note that in Algorithm \ref{alg:mtx} the computation costs of the products $V_M L, G_1 L$, and $G_2 L$, of the QR and SVD decompositions, and of the following truncation step are all independent of $N$.
Secondly, since the nonzero elements of $M_k, D_1, D_2$ are linear in $N$, if Properties P2--P4 above are correct, then
\begin{itemize}
    \item The computation costs of the products $(\sigma_1 \otimes M_k)R, D_1 R$,  and $D_2 R$ scale linearly with $N$;
    \item Computing $K$ is linear in cost with $N$ and the product $RK$ is again linear in $N$ thanks to $K$'s Kronecker structure.
\end{itemize}
Therefore, P1--P4 imply that the cost of Algorithm \ref{alg:mtx} scales linearly with $N$ (note that the discretization cost of the $\star$-strategy is independent of $N$).
Proving Properties P1, P2, and P3 is out of the scope of this paper. However, the numerical evidence we presented supports their correctness.
Table \ref{tab:ctime:a} shows that for a fixed $M$, we achieve the same accuracy for each of the considered $N$ for case (a). Other numerical experiments, not reported here, verify  that this is the case also for (b)--(d) cases. Moreover, the number of singular values above the truncation tolerance is always far below $M$ in all the experiments.
In Section \ref{sec:conv}, we discuss the convergence of the fixed point iteration and its relation to the spectral radius of the iteration matrix. The experiments and the preliminary results on the asymptotic convergence of the method of the next section suggest that the number of iterations is constant or almost constant.
Finally, P4 is a consequence of the other Properties.
\begin{lemma}
    Properties P1, P2, and P3 imply Property P4.
\end{lemma}
\emph{Proof.}
At each iteration, $R$ is a matrix composed of $k \times k$ banded matrices (blocks). In the beginning, $R$ is set to be the $N \times N$ identity matrix that, trivially, is composed of blocks $I_k$ and $k \times k$ null blocks (remember that $N = 2k$).
By induction, given $R$ composed of $k \times k$ banded blocks, we can observe that:
\begin{itemize}
    \item The product $(\sigma_1 \otimes M_k) R$ can cause an increase by at most $1$ of the bandwidth of some of the blocks of $R$, since $\sigma_1 \otimes M_k$ is composed of four $k \times k$ blocks, each of which is either tridiagonal ($M_k$) or a null block. 
    \item The step $R=[D_1 R, D_2 R, d\,]$ does not increase the bandwidth of the blocks as $D_1, D_2$ and $d$ are $N \times N$ diagonal matrices. 
    \item Each of the blocks in $RK$ is obtained by a linear combination of $R$ banded blocks. Therefore, the maximal bandwidth among all blocks of $RK$ is the same as the one of $R$. 
\end{itemize}
Consequently, the maximal bandwidth of each block of $R$ is bounded by the number of iterations, which, by P1, does not asymptotically increase with $N$.

To conclude, note that the size of $R$ is $N \times (rN)$, i.e., the number of $k \times k$ blocks in $R$ is $2r$.
Therefore, as long as $r$ is bounded by $M$ (P3) and $M$ is independent of $N$ (P2), the number of nonzero elements of $R$ is proportional to $N$ times the number of iterations.
$\square$

\subsection{Convergence of the algorithm}\label{sec:conv}
For the stationary iterative method to converge, we noted that the method converges if the spectral radius of the iteration matrix is smaller than 1, i.e., $\rho(G(i\, (\sigma_1\otimes M_k)\otimes V_M))<1$.
That is, the iterates \eqref{eq:it} converge to the solution $x$ of \eqref{eq:sysEq}. 
The asymptotic convergence rate is dictated by the spectral radius \cite{saad03}
\begin{equation*}
    \lim_{m\rightarrow \infty}\left(\frac{\Vert x_m-x\Vert}{\Vert x_0-x\Vert}  \right)^{1/m}= \rho(G(i\, (\sigma_1\otimes M_k)\otimes V_M)).
\end{equation*}
Computing the spectral radius of the matrix is often prohibitive due to its size.
Hence, we will rely on an upper bound that is easier to compute.
It is well known that the spectral radius of a matrix $A$ can be bounded by
\begin{equation*}
    \rho(A)\leq \Vert A^\ell \Vert_2^{1/\ell}.
\end{equation*}
Since multiplication of $G(i\, (\sigma_1\otimes M_k)\otimes V_M)$ with a vector is cheap to compute, we use the Frobenius norm, which can be computed column by column. This leads to a larger upper bound, since
$\Vert A \Vert_2\leq \Vert A\Vert_{\textrm{F}}$, such that
\begin{equation}\label{eq:upperBound}
    \rho(A)\leq \Vert A^\ell \Vert_{\textrm{F}}^{1/\ell}.
\end{equation}
Using a small problem, case (a) of size $N=20$, we illustrate this upper bound and relate the spectral radius to the rate of convergence.
For this small example it is possible to compute the spectral radius.
For $M=130$, we have $\rho(G(i\, (\sigma_1\otimes M_k)\otimes V_M))=0.1780$.
Table \ref{table:upperBound} shows the upper bound \eqref{eq:upperBound} for increasing values of $\ell$.
It shows initially a fast decrease with $\ell$, but this decrease slows down as the upper bound approaches the spectral radius of 0.1780.
\begin{table}[!ht]
		\begin{tabular}{l|llllllll}
			 $\ell$ & $2$  & $4$  & $8$ & $16$  & $32$ & $64$ & $128$ & $256$\\ \hline
			 $\Vert A^\ell\Vert_{\textrm{F}}^{1/\ell}$&  1.97  &  1.16 & 0.806  &  0.578   & 0.415 &   0.296    &0.225  &   0.196 \\
		\end{tabular}
            \caption{Upper bounds on spectral radius of $A:=G(i(\sigma_1\otimes M_k)\otimes V_M)$ for case (a) with $N=20$, which has $\rho(A) = 0.1780$. }
            \label{table:upperBound}
	\end{table}

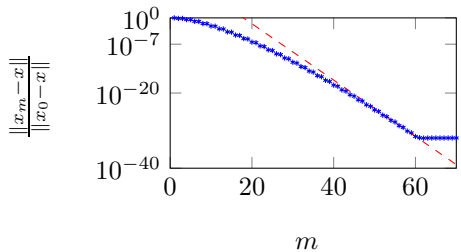
\begin{figure}[!ht]
    \centering
    \setlength\figureheight{2cm}
    \setlength\figurewidth{4cm}
%
\begin{tikzpicture}

\begin{axis}[%
width=0.951\figurewidth,
height=\figureheight,
at={(0\figurewidth,0\figureheight)},
scale only axis,
xmin=0,
xmax=70,
ymode=log,
ymin=1e-40,
ymax=1e+0,
ytick={10^0,10^-7,10^-20,10^-40},
yminorticks=true,
axis background/.style={fill=white}
]
\addplot [color=blue, draw=none, mark=asterisk, mark options={solid, blue}, mark size = 1pt, forget plot]
  table[row sep=crcr]{%
1	1\\
2	0.71\\
3	0.71\\
4	0.35\\
5	0.34\\
6	0.12\\
7	0.11\\
8	0.031\\
9	0.029\\
10	0.0067\\
11	0.0061\\
12	0.0012\\
13	0.0011\\
14	0.00019\\
15	0.00016\\
16	2.6e-05\\
17	2.1e-05\\
18	3.2e-06\\
19	2.4e-06\\
20	3.5e-07\\
21	2.5e-07\\
22	3.5e-08\\
23	2.4e-08\\
24	3.2e-09\\
25	2.1e-09\\
26	2.7e-10\\
27	1.7e-10\\
28	2.1e-11\\
29	1.3e-11\\
30	1.6e-12\\
31	9.2e-13\\
32	1.1e-13\\
33	6.1e-14\\
34	7e-15\\
35	3.8e-15\\
36	4.3e-16\\
37	2.2e-16\\
38	2.5e-17\\
39	1.2e-17\\
40	1.3e-18\\
41	6.4e-19\\
42	7.1e-20\\
43	3.2e-20\\
44	3.5e-21\\
45	1.5e-21\\
46	1.6e-22\\
47	7e-23\\
48	7.5e-24\\
49	3e-24\\
50	3.2e-25\\
51	1.3e-25\\
52	1.3e-26\\
53	5.1e-27\\
54	5.4e-28\\
55	2e-28\\
56	2.1e-29\\
57	7.3e-30\\
58	7.7e-31\\
59	2.6e-31\\
60	2.8e-32\\
61	1.3e-32\\
62	1e-32\\
63	1.1e-32\\
64	1.1e-32\\
65	1.1e-32\\
66	1.1e-32\\
67	1.1e-32\\
68	1.1e-32\\
69	1.1e-32\\
70	1.1e-32\\
};

\addplot [color=red, dashed, forget plot]
  table[row sep=crcr]{%
1	4.8e+12\\
2	8.5e+11\\
3	1.5e+11\\
4	2.7e+10\\
5	4.8e+09\\
6	8.5e+08\\
7	1.5e+08\\
8	2.7e+07\\
9	4.8e+06\\
10	8.6e+05\\
11	1.5e+05\\
12	2.7e+04\\
13	4.8e+03\\
14	8.6e+02\\
15	1.5e+02\\
16	27\\
17	4.9\\
18	0.86\\
19	0.15\\
20	0.027\\
21	0.0049\\
22	0.00087\\
23	0.00015\\
24	2.7e-05\\
25	4.9e-06\\
26	8.7e-07\\
27	1.5e-07\\
28	2.8e-08\\
29	4.9e-09\\
30	8.7e-10\\
31	1.6e-10\\
32	2.8e-11\\
33	4.9e-12\\
34	8.8e-13\\
35	1.6e-13\\
36	2.8e-14\\
37	4.9e-15\\
38	8.8e-16\\
39	1.6e-16\\
40	2.8e-17\\
41	5e-18\\
42	8.8e-19\\
43	1.6e-19\\
44	2.8e-20\\
45	5e-21\\
46	8.9e-22\\
47	1.6e-22\\
48	2.8e-23\\
49	5e-24\\
50	8.9e-25\\
51	1.6e-25\\
52	2.8e-26\\
53	5e-27\\
54	8.9e-28\\
55	1.6e-28\\
56	2.8e-29\\
57	5e-30\\
58	9e-31\\
59	1.6e-31\\
60	2.8e-32\\
61	5.1e-33\\
62	9e-34\\
63	1.6e-34\\
64	2.9e-35\\
65	5.1e-36\\
66	9e-37\\
67	1.6e-37\\
68	2.9e-38\\
69	5.1e-39\\
70	9.1e-40\\
};
\end{axis}

\begin{axis}[%
width=1.227\figurewidth,
height=1.227\figureheight,
at={(-0.16\figurewidth,-0.135\figureheight)},
scale only axis,
xmin=0,
xmax=1,
xlabel= {$m$},
ymin=0,
ymax=1,
ylabel={$\frac{\Vert x_m-x\Vert}{\Vert x_0-x\Vert}$},
axis line style={draw=none},
ticks=none,
axis x line*=bottom,
axis y line*=left,
legend style={legend cell align=left, align=left, draw=white!15!black}
]
\end{axis}
\end{tikzpicture}%
    \caption{Predicted asymptotic convergence $\rho(G(i\, (\sigma_1\otimes M_k)\otimes V_M))^m$ (red dashed line) and observed convergence $\frac{\Vert x_m-x\Vert}{\Vert x_0-x\Vert} $ (blue asterisk) for case (a) with $N=20$.}
    \label{fig:rate_case_a}
 \end{figure}

The observed convergence and theoretically predicted convergence rate, given by the spectral radius, $\rho(G(i\, (\sigma_1\otimes M_k)\otimes V_M))$, are shown in Figure \ref{fig:rate_case_a}.
In this figure, we computed in higher precision, since the observed convergence shows transient behavior until, after 30 iterations, it starts to converge at the theoretically predicted asymptotic rate of convergence.



For larger problems we rely on the Frobenius upper bound.
Table \ref{table:caseA} shows the upper bounds for case (a) for increasing problem size $N$.
For $\ell=16$, the upper bound is below 1, so the fixed point iterations converge.
As $N$ increases, the upper bound increases moderately, notice that the gap between the upper bounds of $N=100$ and $N=500$ is $0.34$ for $\ell=2$ and has decreased to $0.035$ for $\ell=16$.
This suggests that the actual spectral radius does not increase significantly, which implies a similar asymptotic convergence rate.
As noted above, there is transient behavior before the asymptotic convergence rate is attained.
A study of the pseudospectra \cite{TrEm05} would provide more insight into the transient behavior.
In numerical experiments, we observed that the same number of fixed point iterations for all considered $N$ leads to approximately the same accuracy.
\begin{table}[!ht]
\begin{subtable}[t]{0.48\textwidth}
\begin{tabular}{l|lll}
			Case a    &  $\ell=4$  & $\ell=8$  & $\ell=16$ \\ \hline
			$N = 100$ &  1.47 & 0.916 & 0.62 \\
			$N = 200$ & 1.61 & 0.959 & 0.639  \\
			$N = 300$ &  1.69 & 0.985& 0.66 \\
			$N = 400$ & 1.76 & 1.00 & 0.654\\
			$N = 500$ &  1.81 & 1.02 & 0.658
		\end{tabular}
            \caption{}
            \label{table:caseA}
\end{subtable}
\begin{subtable}[t]{0.48\textwidth}
\flushright
\begin{tabular}{l|lll}
			Case b    & $\ell=2$  & $\ell=4$  & $\ell=8$   \\ \hline
			$N = 100$ &2.61   & 1.22    &0.734 \\
			$N = 200$ & 3.13  &  1.33    &0.769  \\
			$N = 300$ &3.47  &  1.40    &0.789 \\
			$N = 400$ & 3.73   & 1.45    &0.804\\
			$N = 500$ &3.94   &  1.50   & 0.815
		\end{tabular}
            \caption{}
            \label{table:caseB}
\end{subtable}

\bigskip 

\begin{subtable}[t]{0.48\textwidth}
\begin{tabular}{l|lll}
			Case c    & $\ell=2$  & $\ell=4$  & $\ell=8$   \\ \hline
			$N = 100$ & 2.63  &  1.21    &0.735 \\
			$N = 200$ &3.14 &   1.33    &0.770 \\
			$N = 300$ &3.48    &1.40     &0.790\\
			$N = 400$ & 3.74 &   1.45   & 0.805\\
			$N = 500$ &3.95 &   1.49    &0.816
		\end{tabular}
            \caption{}
            \label{table:caseC}
\end{subtable}
\hspace{\fill}
\begin{subtable}[t]{0.48\textwidth}
\flushright
\begin{tabular}{l|lll}
			Case d    & $\ell=2$  & $\ell=4$  & $\ell=8$  \\ \hline
			$N = 100$ & 4.77  &  1.09 &   0.514  \\
			$N = 200$ & 5.69  &  1.19  &  0.553 \\
            $N = 300$ &6.30  &  1.25   & 0.55\\
            $N = 400$ &6.78  & 1.30   & 0.563\\
            $N = 500$ &7.18  &  1.34   & 0.571\\    
		\end{tabular}
            \caption{}
            \label{table:caseD}
\end{subtable}
\caption{Upper bounds on spectral radius for case (a)-(d).}
\label{table:allCases}
\end{table}

\subsection{Performance of the cheap error estimator}\label{sec:estimator}
The cost of a numerical method is also determined by our ability to stop it at the right moment. 
We mentioned that our algorithms use a simple and cheap error estimator, here we analyse the performance of this estimator.

\begin{figure}[!ht]
    \centering
    \setlength\figureheight{4cm}
    \setlength\figurewidth{6cm}
%
\begin{tikzpicture}

\begin{axis}[%
width=0.951\figurewidth,
height=\figureheight,
at={(0\figurewidth,0\figureheight)},
scale only axis,
xlabel style={font=\color{white!15!black}},
xlabel={$m$},
ymode=log,
ymin=1e-08,
yminorticks=true,
ylabel style={font=\color{white!15!black}},
ylabel={$\frac{\Vert x_m-x \Vert}{\Vert x_0-x \Vert}$},
axis background/.style={fill=white}
]
\addplot [color=blue, draw=none, mark=asterisk, mark options={solid, blue}, forget plot]
  table[row sep=crcr]{%
1	1.1\\
2	0.67\\
3	0.67\\
4	0.28\\
5	0.27\\
6	0.082\\
7	0.079\\
8	0.021\\
9	0.019\\
10	0.004\\
11	0.0037\\
12	0.00074\\
13	0.00063\\
14	0.00011\\
15	9.1e-05\\
16	1.5e-05\\
17	1.2e-05\\
18	1.8e-06\\
19	1.4e-06\\
20	2e-07\\
21	1.5e-07\\
};
\addplot [color=red, draw=none, mark=o, mark options={solid, red}, forget plot]
  table[row sep=crcr]{%
1	89\\
2	20\\
3	0.22\\
4	5.3\\
5	0.07\\
6	1.2\\
7	0.018\\
8	0.23\\
9	0.0036\\
10	0.039\\
11	0.00064\\
12	0.0058\\
13	9.7e-05\\
14	0.00077\\
15	1.3e-05\\
16	9.1e-05\\
17	1.6e-06\\
18	9.8e-06\\
19	1.7e-07\\
20	9.7e-07\\
21	1.7e-08\\
};

\end{axis}
\end{tikzpicture}%
    \caption{Observed convergence (blue asterisk) and residual estimator (red circle) for case (a) with $N=400$.}
    \label{fig:convRate:case_a:d=400}
 \end{figure}
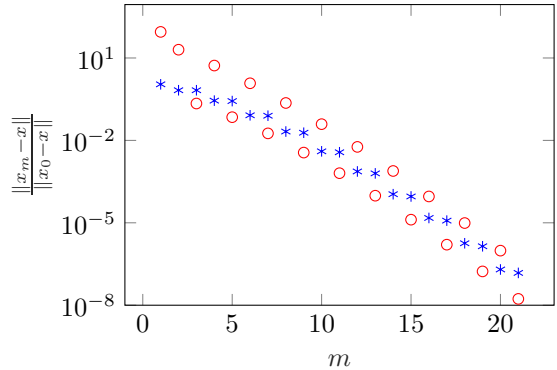

In Figure \ref{fig:convRate:case_a:d=400} the error during the fixed point iterations is shown together with the cheap error estimator used for our algorithms.
This cheap error estimator describes the actual error quite accurately.
Note that an error of about $10^{-7}$ is obtained after $m=20$ iterations, which we also observed for the smaller problem $N=20$ in Figure \ref{fig:rate_case_a}.

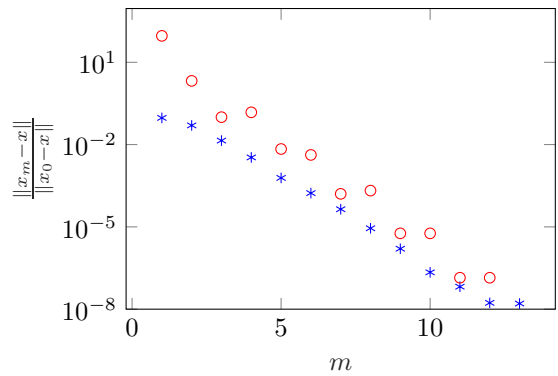
\begin{figure}[!ht]
    \centering
    \setlength\figureheight{4cm}
    \setlength\figurewidth{6cm}
%
\begin{tikzpicture}

\begin{axis}[%
width=0.951\figurewidth,
height=\figureheight,
at={(0\figurewidth,0\figureheight)},
scale only axis,
xlabel style={font=\color{white!15!black}},
xlabel={$m$},
ymode=log,
ymin=1e-08,
yminorticks=true,
ylabel style={font=\color{white!15!black}},
ylabel={$\frac{\Vert x_m-x \Vert}{\Vert x_0-x \Vert}$},
axis background/.style={fill=white}
]
\addplot [color=blue, draw=none, mark=asterisk, mark options={solid, blue}, forget plot]
  table[row sep=crcr]{%
1	0.094\\
2	0.05\\
3	0.014\\
4	0.0034\\
5	0.00061\\
6	0.00017\\
7	4.4e-05\\
8	8.9e-06\\
9	1.6e-06\\
10	2.2e-07\\
11	6.6e-08\\
12	1.7e-08\\
13	1.6e-08\\
};
\addplot [color=red, draw=none, mark=o, mark options={solid, red}, forget plot]
  table[row sep=crcr]{%
1	92\\
2	2.1\\
3	0.1\\
4	0.15\\
5	0.0069\\
6	0.0042\\
7	0.00016\\
8	0.00021\\
9	5.8e-06\\
10	5.8e-06\\
11	1.4e-07\\
12	1.4e-07\\
13	2.6e-09\\
};

\end{axis}

\begin{axis}[%
width=1.227\figurewidth,
height=1.227\figureheight,
at={(-0.16\figurewidth,-0.135\figureheight)},
scale only axis,
xmin=0,
xmax=1,
ymin=0,
ymax=1,
axis line style={draw=none},
ticks=none,
axis x line*=bottom,
axis y line*=left,
legend style={legend cell align=left, align=left, draw=white!15!black}
]
\end{axis}
\end{tikzpicture}%
    \caption{Observed convergence (blue asterisk) and residual estimator (red circle) for case (d) with $N=400$.}
    \label{fig:convRate:case_d:d=400}
 \end{figure}

For cases (b), (c) and (d) the cheap error estimator also describes the actual error sufficiently well.
Case (d) is shown in Figure \ref{fig:convRate:case_d:d=400}.

\section{Conclusions}\label{sec:conclusions}
In this paper, we have introduced a new approach for the solution of the Rosen-Zener non-autonomous linear system of ODEs. 
The new approach is based on the so-called $\star$-product expression for the solution of ODEs and on the related discretization obtained by Legendre orthogonal polynomials.
In the more challenging problem of computing the operator solution of the equation, the new Algorithm \ref{alg:mtx} appears to scale linearly with the size of the system. The experiments showed that in general for large systems and for solutions on large intervals the new approach is much faster than state-of-the-art methods. However, when high accuracy (close to machine precision) is needed, or when the application requires preserving geometrical properties (e.g., the unitarity of the solution), then our approach might fail, and the geometrical integration methods we used for comparison can be a better option. In particular, at present, geometric properties seem not to be preserved by our approach. 
Since our proposed method is a global method, i.e., it computes a single approximating polynomial on the whole time-domain of interest, its deviation from unitarity is expected to be of the order of the error of the approximation.
This error can be estimated by looking at the amplitude of the Legendre coefficients representing this polynomial.
As opposed to time-stepping methods there is no build up of error by summing many local approximations on subintervals. For time-stepping methods, if in every step there is a small deviation from unitarity, the approximation at the final time might deviate significantly from unitarity.
Hence, for our global method, solving the same problem on a larger time interval will not necessarily increase the deviation from unitarity, on the condition that the condition number of the discretized matrix does not increase significantly and the available computer resources allow us to solve the problem on the whole interval of interest.
Moreover, the approximation of a smooth function by a single polynomial is known to exhibit spectral convergence, which is the fastest rate of convergence possible and cannot be achieved by the piecewise approximations implicit in time stepping methods.

One of the key points is that the new $\star$-approach allows for the exploitation of hidden structural properties of the problem, such as the low numerical rank of the related matrix equation's solution. These properties made such a fast algorithm possible. The full connection between the matrix equation properties and the cost of the algorithm requires proving Properties~P1--P3 in Section \ref{sec:cost}. While substantial numerical evidence for these Properties has been provided in the paper, a complete numerical analysis of the algorithm will be developed in future work.

Algorithm \ref{alg:mtx} exploits structure particular to the Rosen-Zener model. However, the underlying approach can be adapted to other, more general models, once a fast solver for the related matrix equation is identified. We are currently working on these generalizations.

\section*{Appendix}
\begin{table}[h!]
\centering
\begin{tabular}{ |c||c|c|c|c|c|c|c|  }
 \hline
   & star & $\text{SM}_8^{[4]}$ & $\text{SM}_{11}^{[6]}$ & $\text{SM}_{11}^{[8]}$ & $\text{S}_6^{[4]}$ & $\text{S}_{10}^{[6]}$ & $\text{RK}_7^{[6]} $ \\
 \hline
 \hline
  Case   &  \multicolumn{7}{|c|}{Method's parameters} \\
    \hline
  (a)   & $M=130$ & $\kappa = 8$ & $\kappa = 7$ & $\kappa = 6$ & $\kappa = 10$ & $\kappa = 8$ & $\kappa = 10$ \\
  \hline
  \hline
  $N$ &  \multicolumn{7}{|c|}{Norm of the absolute error $\times 1e-7$} \\
 \hline 
 160 & 1.506 &   8.789 &    0.509 &    1.846 &    4.112 &    2.721 &    2.970 \\
 320 & 1.525 &   8.789 &    0.509 &    1.846 &    4.112 &    2.721 &    2.970 \\
 480 & 1.529 &   8.789 &    0.509 &    1.846 &    4.112 &    2.721 &    2.970 \\
 640 & 1.530 &   8.789 &    0.509 &    1.846 &    4.112 &    2.721 &    2.970 \\
 800 & 1.531 &   8.789 &    0.509 &    1.846 &    4.112 &    2.721 &    2.970 \\
 $\vdots$ & $\vdots$ & $\vdots$ & $\vdots$ & $\vdots$ & $\vdots$ & $\vdots$ & $\vdots$ \\ 
1600 & 1.531 &   8.789 &    0.509 &    1.846 &    4.112 &    2.721 &    2.970 \\ 
 \hline
\end{tabular}
\caption{Experiment 2, case (a). Parameter settings of numerical methods and maximal absolute error on the interval for increasing $N$. For all methods, except star, the interval is split into $2^\kappa$ subintervals.}\label{tab:ctime:a}
\end{table}
\begin{table}[h!]
\centering
\begin{tabular}{ |c||c|c|c|c|c|c|c|  }
 \hline
   & star & $\text{SM}_8^{[4]}$ & $\text{SM}_{11}^{[6]}$ & $\text{SM}_{11}^{[8]}$ & $\text{S}_6^{[4]}$ & $\text{S}_{10}^{[6]}$ & $\text{RK}_7^{[6]} $ \\
 \hline
 \hline
  Cases  &  \multicolumn{7}{|c|}{Methods' parameters} \\
    \hline
  (b)   & $M=130$ & $\kappa = 8$ & $\kappa = 7$ & $\kappa = 6$ & $\kappa = 10$ & $\kappa = 8$ & $\kappa = 10$ \\
  (c)   & $M=210$ & $\kappa = 9$ & $\kappa = 7$ & $\kappa = 7$ & $\kappa = 10$ & $\kappa = 8$ & $\kappa = 10$ \\
  (d)   & $M=500$ & $\kappa = 10$ & $\kappa = 9$ & $\kappa = 8$ & $\kappa = 10$ & $\kappa = 9$ & $\kappa = 10$ \\
  \hline
  \hline
   &  \multicolumn{7}{|c|}{Norm of the absolute error $\times 1e-7$, $N=1600$} \\
 \hline 
 (b) & 0.874 &   8.909 &    0.553 &    6.530 &    4.109 &    3.038 &    3.057 \\
 (c) & 0.808 &   0.725 &    0.631 &    2.416 &    4.215 &    2.958 &    3.199 \\
 (d) & 0.156 &   6.706 &    0.172 &    6.511 &    3.424 &    0.247 &    8.142 \\
 \hline
\end{tabular}
\caption{Experiment 2, case (b)--(d). Parameter settings of numerical methods and maximal absolute error on the interval for $N=1600$. For all methods, except star, the interval is split into $2^\kappa$ subintervals..}\label{tab:ctime:bcd}
\end{table}


\begin{table}[h!]
\centering
\begin{tabular}{ |c||c|c|c|c|c|c|c|  }
 \hline
   & star & $\text{SM}_8^{[4]}$ & $\text{SM}_{11}^{[6]}$ & $\text{SM}_{11}^{[8]}$ & $\text{S}_6^{[4]}$ & $\text{S}_{10}^{[6]}$ & $\text{RK}_7^{[6]} $ \\
 \hline
 \hline
 Case (a)&  \multicolumn{7}{|c|}{Methods' parameters} \\
    \hline
  $I$ length & $M$ & $\kappa$ & $\kappa$ & $\kappa$ & $\kappa$ & $\kappa$ & $\kappa$   \\
    \hline
   $50.2$  & $210$ & $10$ & $8$ & $7$ & $11$ & $9$ & $11$ \\
   $100.5$ & $370$ & $11$ & $9$ & $8$ & $13$ & $11$ & $13$  \\
   $150.7$  & $530$ & $12$ & $9$ & $9$ & $13$ & $11$ & $13$ \\
   $201.0$  & $690$ & $12$ & $10$ & $9$ & $14$ & $12$ & $14$  \\
   $251.3$  & $850$ & $13$ & $10$ & $10$ & $14$ & $12$ & $14$  \\
  \hline
   &  \multicolumn{7}{|c|}{Norm of the absolute error $\times 1e-7$} \\
   \hline
   $50.2$  &     1.794  &  1.103 &   0.528  &  1.974  &  7.911  &  5.400  &  5.753 \\
   $100.5$ & 0.178  &  2.172 &   0.563  &  1.908  &  0.970  &  0.166  &  0.177 \\
   $150.7$  & 0.178  &  1.031 &   6.734  &  0.244  &  7.315  &  2.854  &  3.004 \\
   $201.0$  & 0.178  &  4.312 &   0.652  &  1.761  &  1.920  &  0.332  &  0.351 \\
   $251.3$  & 0.178  &  0.826 &   2.620  &  0.071  &  5.849  &  1.588  &  1.669 \\  
 \hline 
\end{tabular}
\caption{Parameter settings of numerical methods for Experiment 4, Case (a), and the absolute error of the operator approximation for different lengths of the interval $I$. For all methods, except star, the interval is split into $2^\kappa$ subintervals.}\label{tab:interval:a}
\end{table}   
\begin{table}[h!]
\centering
\begin{tabular}{ |c||c|c|c|c|c|c|c|  }
 \hline
   & star & $\text{SM}_8^{[4]}$ & $\text{SM}_{11}^{[6]}$ & $\text{SM}_{11}^{[8]}$ & $\text{S}_6^{[4]}$ & $\text{S}_{10}^{[6]}$ & $\text{RK}_7^{[6]} $ \\
 \hline
 \hline
   Case (d) &  \multicolumn{7}{|c|}{Methods' parameters} \\
   \hline
  $I$ length & $M$ & $\kappa$ & $\kappa$ & $\kappa$ & $\kappa$ & $\kappa$ & $\kappa$   \\
    \hline
   $12.5$  & $230$ & $9$ & $7$ & $7$ & $9$ & $7$ & $9$ \\
   $25.1$ & $410$ & $10$ & $9$ & $8$ & $10$ & $9$ & $10$  \\
   $37.6$  & $600$ & $11$ & $9$ & $9$ & $11$ & $9$ & $11$ \\
   $50.2$  & $800$ & $12$ & $10$ & $10$ & $11$ & $10$ & $12$  \\
   $62.8$  & $1000$ & $12$ & $10$ & $10$ & $12$ & $10$ & $12$  \\
   $75.3$  & $1200$ & $12$ & $10$ & $10$ & $12$ & $10$ & $12$  \\
   $87.9$  & $1400$ & $13$ & $11$ & $10$ & $12$ & $11$ & $13$  \\
 \hline
   &  \multicolumn{7}{|c|}{Norm of the absolute error $\times 1e-7$ } \\
   \hline
   $12.5$  & 0.641  &  3.353  &  5.472  &  3.182  &  1.712  &  7.915  &  4.087 \\
   $25.1$ & 0.330  &  6.706  &  0.172  &  6.510  &  3.424  &  0.247  &  8.142 \\
   $37.6$  & 0.977  &  3.189  &  2.932  &  1.749  &  1.624  &  4.226  &  2.173 \\
   $50.2$  & 0.101  &  0.841  &  0.344  &  0.205  &  6.848  &  0.495  &  0.254 \\
   $62.8$  & 0.142  &  2.565  &  1.639  &  0.983  &  1.305  &  2.358  &  1.212 \\
   $75.3$  & 0.146  &  6.378  &  5.863  &  3.531  &  3.249  &  8.451  &  4.342 \\
   $87.9$  & 0.152  &  0.863  &  0.270  &  0.162  &  7.023  &  0.388  &  0.200 \\
 \hline 
\end{tabular}
\caption{Parameter settings of numerical methods for Experiment 4, Case (d), and the absolute error of the operator approximation for different lengths of the interval $I$. For all methods, except star, the interval is split into $2^\kappa$ subintervals.}\label{tab:interval:d}
\end{table}  

\newpage


\end{document}